\documentclass[pdflatex,sn-mathphys-num]{sn-jnl}

\usepackage{graphicx}%
\usepackage{multirow}%
\usepackage{amsmath,amssymb,amsfonts}%
\usepackage{amsthm}%
\usepackage{mathrsfs}%
\usepackage[title]{appendix}%
\usepackage{xcolor}%
\usepackage{textcomp}%
\usepackage{manyfoot}%
\usepackage{booktabs}%
\usepackage{algorithm}%
\usepackage{algorithmicx}%
\usepackage{algpseudocode}%
\usepackage{listings}%

\theoremstyle{thmstyleone}%
\newtheorem{theorem}{Theorem}
\newtheorem{corollary}[theorem]{Corollary}%
\newtheorem{lemma}[theorem]{Lemma}%

\theoremstyle{thmstyletwo}%
\newtheorem{remark}{Remark}%

\theoremstyle{thmstylethree}%

\newcommand{\C}{\mathbb{C}}
\newcommand{\D}{\mathbb{D}}
\newcommand{\R}{\mathbb{R}}
\newcommand{\N}{\mathbb{N}}
\renewcommand{\Im}{\mathop{\textnormal{Im}}}

\DeclareMathOperator{\id}{\textnormal{id}}
\DeclareMathOperator{\Hol}{\textnormal{Hol}}
\DeclareMathOperator{\Aut}{\textnormal{Aut}}

\begin{document}

\title[Parabolic dynamics according to Maurice Heins]{Parabolic dynamics according to Maurice Heins}


\author*[1,2]{\fnm{Marco} \sur{Abate}}\email{marco.abate@unipi.it}

%

\affil*[1]{\orgdiv{Dipartimento di Matematica}, \orgname{Universit\`a di Pisa}, \orgaddress{\street{Largo Pontecorvo 5}, \city{Pisa}, \postcode{56127}, \state{(PI)}, \country{Italy}}}

\affil[2]{\orgdiv{} \orgname{Universit\`a degli Studi Guglielmo Marconi}, \orgaddress{\street{Via Plinio 44}, \city{Roma}, \postcode{00193}, \state{(Roma)}, \country{Italy}}}



\abstract{This paper, based upon an unpublished manuscript by Maurice Heins, answers a question posed by Valiron about the dynamics of parabolic self-maps 
of the unit disk in the complex plane, considerably simplifying arguments previously used for answering the same question. The main new tool introduced is the notion of left straightening of a sequence of iterates, that can be effectively employed for studying the hyperbolic step of a parabolic map.}




\maketitle

\vspace*{-1cm}
\centerline{\textit{This paper is dedicated to the memory of Marek Jarnicki}}
\bigskip

\bmhead{Acknowledgements}
The author was partially supported by the French-Italian University and by Campus France through the Galileo program, under the project \textit{From rational to transcendental: complex dynamics and parameter spaces}, as well as by Istituto Nazionale di Alta Matematica (INdAM).
\bigskip

\section{Introduction}\label{section1}

Marek Jarnicki was a great expert on invariant distances and metrics in complex analysis and geometry; his book with Peter Pflug \cite{JarPfl2013} is a definitive reference on the subject. Since the 1920's, invariant distances have often been instrumental in the study of holomorphic dynamics on hyperbolic Riemann surfaces and hyperbolic manifolds; a good example is Wolff's use \cite{Wol1926} of the result that has later become known as the Wolff lemma to simplify the proof of the Wolff-Denjoy theorem (see Theorem~\ref{th:WD}) describing the dynamics of holomorphic self-maps of the unit disk~$\mathbb{D}$ of the complex plane. 

Maurice Heins (1915--2015) was a master in the study of holomorphic dynamics on hyperbolic Riemann surfaces; for instance, he has been able \cite{He1941} to extend the Wolff-Denjoy theorem to multiply connected domains in the plane (for more on this beautiful theory and its history see, e.g., \cite{Ab2022} and references therein).

He often used invariant distances in his work, as confirmed by a recent unexpected (at least by me) discovery. 
While Ian Short and myself were working on the paper~\cite{AbSh2025} on iterated function systems (where the object of study is the dynamics of a sequence of functions obtained by composing different self-maps), we came across (actually, Ian did) an unpublished manuscript by Heins \cite{He1991}, where with a clever use of the Poincar\'e invariant distance he proved a number of beautiful results on holomorphic dynamics in~$\D$. We incorporated, extended and generalized some of these results in~\cite{AbSh2025}, but a few of them were left out, because they were not relevant to our setting. In particular, we left out an elegant solution to a question on parabolic dynamics posed by Valiron in~\cite{Val1954}. This question has already been answered by Pommerenke in~\cite{Pomm1979}, but Heins' approach is surprisingly simpler and based on ideas having other applications too (see, for instance, \cite{AbSh2025} and Section~\ref{section2} of this paper). 

I believe that Heins' ideas deserve to be better known; being based on invariant distances, they fit well in a publication dedicated to Marek. So this paper is a report 
on Heins' approach to parabolic dynamics as described in~\cite{He1991}. I have updated the presentation, in order to make it coherent with contemporary research on this topic, and streamlined some proofs; but the main ideas are his.

Let me now summarize what we are going to discuss in this paper. Let $\mathbb{D}=\{z\in\C\mid |z|<1\}$ be the unit disk in the complex plane.
We denote by $\Hol(\D,\D)$ the space of holomorphic self-maps of~$\D$, by $\Aut(\D)$ the set of (holomorphic) automorphisms of~$\D$ and by~$\omega$ the Poincar\'e distance of~$\D$. Finally, given $f\in\Hol(\D,\D)$, we denote by $f^n$ its $n$-th iterate, the composition of~$f$ with itself $n$ times. Finally, the orbit of a point $z\in\D$ is the sequence~$\{f^n(z)\}$. 

If $f\in\Hol(\D,\D)$ has a fixed point $z_0\in\D$, an easy application of the Schwarz-Pick lemma shows that either $f$ is an automorphism (and  then its dynamics is the same as the dynamics of an Euclidean rotation) or the sequence of iterates converges to~$z_0$. Furthermore, thanks to the work of K\oe nigs \cite{Koe} in 1883 and B\"ottcher \cite{Bot} in 1904, we know exactly how the orbits converge to the fixed point.

The fundamental result describing the dynamics of fixed point free self-maps is the Wolff-Denjoy theorem mentioned above, that can be summarized as follows:

\begin{theorem}
\label{th:WD}
Let $f\in\Hol(\D,\D)$ be a holomorphic self-map without fixed points. Then there exists a point $\tau_f\in\partial\D$ such that $f^n\to\tau_f$ uniformly on compact subsets. Furthermore the derivative $f'$ admits non-tangential limit $f'(\tau_f)$ at~$\tau_f$ and $f'(\tau_f)\in(0,1]$.
\end{theorem}

The point $\tau_f$ is the \emph{Wolff point} of the map $f$. 
Using the position of the Wolff point and the value of $f'(\tau_f)$ we can introduce the following classification: a holomorphic map $f\in\Hol(\D,\D)$ is
\begin{enumerate}
\item \emph{elliptic} if it has a fixed point in~$\D$;
\item \emph{parabolic} if it has no fixed points in $\D$ and $f'(\tau_f)=1$;
\item \emph{hyperbolic} if it has no fixed points in $\D$ and $0<f'(\tau_f)<1$.
\end{enumerate} 
The dynamics of hyperbolic maps is more complicated than the dynamics of elliptic maps but still well-behaved: Wolff \cite{Wol1929} in 1929 and Valiron \cite{Val1931} in 1931 proved that for any $z_0\in\D$ the orbit $\{f^n(z_0)\}$ converges to the Wolff point non-tangentially and with a precise slope. More precisely, Valiron~\cite{Val1931} (see also \cite{Val1954}) proved that if $z_0\in\D$ is fixed then there exists a non-constant holomorphic map~$\psi$ such that
\begin{equation}
\lim_{n\to+\infty}\frac{f^n(z)-\tau_f}{f^n(z_0)-\tau_f}=\psi(z)
\label{eq:int0}
\end{equation}
and, moreover, $\psi$ is a solution of the Schr\"oder equation $\psi\circ f=f'(\tau_f)\psi$. Furthermore, Valiron also proved that
for every $z\in\D$ there exists $\theta_z\in(-\pi/2,\pi/2)$ such that
\begin{equation}
\lim_{n\to+\infty}\frac{f^n(z)-\tau_f}{|f^n(z)-\tau_f|}=\tau_f e^{i\theta_z}\;.
\label{eq:int1}
\end{equation}
In the parabolic case, Wolff and Valiron were able to obtain partial results only. So in \cite[p. 148]{Val1954} Valiron asked whether the limit of the left-hand side in \eqref{eq:int0}
does exist in the parabolic case too. Moreover, he wondered whether there are conditions ensuring that parabolic orbits converge to the Wolff point along a definite slope, like in~\eqref{eq:int1}. (Actually, Valiron worked in the right half-plane; the previous statements are the translation of Valiron's results and questions to~$\D$, which is biholomorphic to a half-plane in~$\C$ via the Cayley transform; see Section~\ref{section3}.) 

As mentioned before, Pommerenke \cite{Pomm1979} answered these questions in~1979. However, his answers are a byproduct of quite a delicate argument whose aim is to find a solution to the Abel equation $\psi\circ f=\psi+b$ or, in modern terminology, to find models for parabolic maps. Here, a \emph{model} is an automorphism of~$\D$ or~$\C$ which is semiconjugate in an appropriate sense to the original map~$f$ and thus can be used to study its dynamics. The theory of models is very powerful but definitely not easy; see \cite[Chapters~3 and 4]{Ab2022} for an introduction.

It is worthwhile to remark immediately that Wolff and Valiron already knew that \eqref{eq:int1} does \textit{not} hold for all parabolic maps. More precisely, Pommerenke clarified that the class of parabolic maps should be subdivided in two subclasses, having quite a different dynamical behaviour. In modern terminology, the subdivision is expressed in terms of the hyperbolic step, that, in turn, is defined by using the Poincar\'e distance~$\omega$ on~$\D$. Take $f\in\Hol(\D,\D)$. By the Schwarz-Pick theorem, for any $z\in\D$ the sequence $\bigl\{\omega\bigl(f^n(z),f^{n+1}(z)\bigr)\bigr\}$ is non-increasing and, thus, it is converging. The \emph{hyperbolic step} $s^f\colon\D\to\R^+$ of~$f$ is then defined by
\[
s^f(z)=\lim_{n\to+\infty}\omega\bigl(f^n(z),f^{n+1}(z)\bigr)\;.
\]
It turns out that the dynamics of parabolic maps with positive hyperbolic step (that is, such that $s^f>0$ everywhere) is quite different from the dynamics of parabolic maps with zero hyperbolic step (that is, such that $s^f\equiv 0$). In particular, \eqref{eq:int1} holds for parabolic maps with positive hyperbolic step but there are examples of parabolic maps with zero parabolic step where \eqref{eq:int1} does not hold. 

In principle, there might be a third class of parabolic maps, consisting of maps for which the hyperbolic step vanishes somewhere but it is not identically zero. However, it turns out that this is not the case: $s^f$ vanishes at one point if and only if it is identically zero. This result, which is instrumental in proving \eqref{eq:int1} for parabolic maps with positive hyperbolic step, has been obtained by using the theory of models; however, in~\cite[p. 248]{Ab2022} I wondered whether it was possible to give a proof independent of the theory of models.

Heins' answers to these questions is based on the apparently unrelated idea of left straightening of the sequence of iterates~$\{f^n\}$ of a map~$f$ (here I am using the terminology introduced, in a more general setting, in~\cite{AbSh2025}). A \emph{left straightening} of $\{f^n\}$ is a holomorphic self-map $h\in\Hol(\D,\D)$ obtained as limit of a sequence of the form $\{\gamma_n^{-1}\circ f^n\}$, with $\gamma_n\in\Aut(\D)$ for all $n\in\N$. 

In Section~\ref{section2} we shall prove that a left straightening always exists and is unique up to left composition by an automorphism (Theorem~\ref{th:leftst}). 
Since the automorphisms of~$\D$ are isometries for the Poincar\'e distance, a left straightening~$h$ of $\{f^n\}$ can be used to express the hyperbolic step: indeed, we have
$s^f(z)=\omega\bigl(h(z),h\bigl(f(z)\bigr)$ for all $z\in\D$ (Corollary~\ref{th:lshs}). So $s^f(z_0)=0$ if and only if $h\bigl(f(z_0)\bigr)=h(z_0)$. Having this, a clever application of Rouch\'e's theorem allows to prove that $s^f$ vanishes at one point if and only if it is identically zero (Theorem~\ref{th:zerohypst}).  

Finally, in Section~\ref{section3} we shall show how to use the hyperbolic step and Theorem~\ref{th:zerohypst} to push the original Valiron's argument just that little bit further needed to answer his questions: when $f$ is parabolic, then the ratio in \eqref{eq:int0} has limit~1 (Corollary~\ref{th:Val1D}) and when $p$ is parabolic with positive hyperbolic step then the ratio in~\eqref{eq:int1} has the same limit ~$\pm i\tau_f$ for all~$z\in\D$ (Corollary~\ref{th:Val2D}). The final argument is deceptively simple; but this is just a consequence of the elegance and strength of Heins' ideas.

\section{Left straightening and hyperbolic step}
\label{section2}

Let us fix a few notations and recall a few well-know results; the details can be found in \cite{Ab2022} and references therein.

The Poincar\'e distance $\omega\colon\D\times\D\to\R^+$ is given by
\[
\omega(z_1,z_2)=\tanh^{-1}\left|\frac{z_2-z_1}{1-\overline{z_1}z_2}\right|\;,
\]
where $\tanh^{-1} t=\frac{1}{2}\log\frac{1+t}{1-t}$; in particular, $\omega(0,z)=\tanh^{-1}|z|$.

%
%
The fundamental Schwarz-Pick theorem says that holomorphic maps are semicontractions with respect to the Poincar\'e distance: if $f\in\Hol(\D,\D)$ then
\[
\omega\bigl(f(z),f(w)\bigr)\le\omega (z,w)
\]
for all $z$,~$w\in\D$, with equality for some $z\ne w$ if and only if equality holds for all $z$,~$w\in\D$ if and only if $f\in\Aut(\D)$. 

As mentioned in the introduction, a consequence of the Schwarz-Pick theorem is that, for any $z\in\D$, the sequence $\bigl\{\omega\bigl(f^n(z),f^{n+1}(z)\bigr)\bigr\}$ is non-increasing and, thus, it is converging. This allows us to define the \emph{hyperbolic step} $s^f\colon\D\to\R^+$ by
\[
s^f(z)=\lim_{n\to+\infty}\omega\bigl(f^n(z),f^{n+1}(z)\bigr)\;.
\]
We shall say that $f$ has \emph{positive hyperbolic step} if there exists $z_0\in\D$ such that $s^f(z_0)>0$; and that $f$ has \emph{zero hyperbolic step} if instead $s^f\equiv 0$. 

One natural question now is whether there might exist self-maps (not automorphisms; see below) whose hyperbolic step is vanishing somewhere but it is not identically zero. The answer is negative. This has been remarked by many authors, including Pommerenke; but their proofs depended (implicitly or explicitly) on elaborated computations underlying the theory of models for holomorphic self-maps. 

In this section we shall give a simple proof of this fact, by using the idea of left straightening of a sequence of holomorphic self-maps, implicitly introduced by Heins in the unpublished manuscript~\cite{He1991} for the sequence of iterates of a single map and then generalized by Short and the author~\cite{AbSh2025} to arbitrary iterated function systems of~$\D$. For the sake of completeness, we report here a proof of the existence of a left straightening of a sequence of iterates (for the general case see \cite[Theorem A]{AbSh2025}). The original proof by Heins was based on the Harnack convergence theorem; here we instead use a normal families argument, which is more in line with the rest of the paper.

\begin{theorem}
\label{th:leftst}
Given $f\in\Hol(\D,\D)$, there exists a holomorphic map $h\in\Hol(\D,\D)$ and a sequence $\{\gamma_n\}\subset\Aut(\D)$ such that $\gamma_n^{-1}\circ f^n$ converges to~$h$ uniformly on compact subsets. The map $h$ is unique up to left composition by an automorphism of~$\D$. Moreover, given $z_0\in\D$ we can choose $\{\gamma_n\}$ so that $\gamma_n^{-1}\circ f^n(z_0)=h(z_0)=0$ for all~$n\in\N$.
\end{theorem}

\begin{proof}
Given $z_0\in\D$, we choose $\gamma_n\in\Aut(\D)$ with $\gamma_n(0)=f^n(z_0)$ and let $H_n=\gamma_n^{-1}\circ f^n$. Then $H_n(z_0)=0$; in particular, $\{H_n\}$ is a relatively compact family in~$\Hol(\D,\D)$. For $z\in\D$, we have
\[
\omega\bigl(H_n(z),0\bigr) = \omega\bigl(H_n(z),H_n(z_0)\bigr)=\omega\bigl(f^n(z),f^n(z_0)\bigr)\;.
\]
Since $\omega\bigl(f^n(z),f^n(0)\bigr)\le \omega\bigl(f^{n-1}(z),f^{n-1}(0))$, it follows that $\omega\bigl(H_n(z),0\bigr)\le \omega\bigl(H_{n-1}(z),0\bigr)$; therefore $\{|H_n|\}$  
is a non-increasing sequence. If $|H_n|\to 0$ (pointwise and hence, by Vitali's theorem, uniformly on compact subsets), then we can take $h\equiv 0$ and we are done.
Otherwise, there exists $w_0\in \D$ for which $\{|H_n(w_0)|\}$ converges to a positive constant. For any $n\in\N$, let $\theta_n$ be an argument of $H_n(w_0)$. By pre-composing $\gamma_n$ with the rotation $z\mapsto e^{i\theta_n}z$, we can assume that $\{H_n(w_0)\}$ is a non-increasing sequence of positive numbers converging to a positive number~$\rho_0\in(0,1)$.

Suppose now that there are two subsequences, $\{H_{m_i}\}$ and $\{H_{n_j}\}$, of $\{H_n\}$ converging to $h$ and $k$ respectively. We clearly have $h(z_0)=k(z_0)=0$ and $h(w_0)=k(w_0)=\rho_0$. By passing to further subsequences we can assume that $m_1<n_1<m_2<n_2<\dotsb$. Let $K_i=\gamma_{n_i}^{-1}\circ f^{n_i-m_i}\circ \gamma_{m_i}$. Then 
\begin{equation}
H_{n_i}=K_i\circ H_{m_i}\;.
\label{eq:ast}
\end{equation} 
Note that $K_i(z_0)=0$; so $\{K_i\}$ is relatively compact in $\Hol(\D,\D)$. Consequently, there is a subsequence of~$\{K_i\}$ converging to $\psi\in \Hol(\mathbb{D},\mathbb{D})$; from \eqref{eq:ast} we infer that $\psi \circ h=k$. Notice that 
\[
\psi(0)=\psi\bigl(h(z_0)\bigr)=k(z_0)=0\;;
\] 
analogously, $\psi(\rho_0)=\rho_0$. It follows that $\psi=\id_\D$ since, among all holomorphic self-maps of~$\mathbb{D}$, only the identity map fixes two distinct points (this is a well-known easy consequence of the uniqueness part of the Schwarz-Pick lemma); hence, $h\equiv k$. 
Thus the relatively compact sequence $\{H_n\}$ has a unique limit point~$h$ and, thus, $H_n\to h$, as claimed.

It remains to prove that $h$ is unique up to left composition by elements of $\Aut(\D)$. Suppose then that there are sequences $\{\gamma_n\}$ and $\{\delta_n\}$ in $\Aut(\D)$, such that $\gamma_n^{-1}\circ f^n\to h$ and $\delta_n^{-1}\circ f^n\to k$, with $h$,~$k\in\Hol(\D,\D)$. Let $\phi_n=\gamma_n^{-1}\circ \delta_n$; here we are not assuming anything on the value of $h$ and $k$ in~$z_0$. Then
\[
\begin{aligned}
\omega\bigl(\phi_n\bigl(k(z_0)\bigr),h(z_0)\bigr) 
&\le \omega\bigl(\phi_n\bigl(k(z_0)\bigr),\phi_n\bigl(\delta_n^{-1}\bigl(f^n(z_0)\bigr)\bigr)\bigr)+\omega\bigl(\phi_n\bigl(\delta_n^{-1}\bigl(f^n(z_0)\bigr)\bigr),h(z_0)\bigr)\\
& = \omega\bigl(k(z_0),\delta_n^{-1}\circ f^n(z_0)\bigr)+\omega\bigl(\gamma_n^{-1}\circ f^n(z_0),h(z_0)\bigr)\;.
\end{aligned}
\]
Hence $\omega\bigl(\phi_n\bigl(k(z_0)\bigr),h(z_0)\bigr)\to 0$ and, thus, $\{\phi_n\}$ is relatively compact in~$\Aut(\D)$. It follows that it has a subsequence converging to $\phi\in \Aut(\D)$. From $\gamma_n^{-1}\circ f^n = \phi_n\circ (\delta_n^{-1}\circ f^n)$, passing to the limit along this subsequence we obtain $h=\phi\circ k$, as required. 
\end{proof}

A map $h$ as given by the previous statement is a \emph{left straightening} of the sequence $\{f^n\}$ of iterates of~$f$. 

In \cite{AbSh2025} we used the left straightening to prove several results on iterated function systems. Here, the main point is that we can use any left straightening to compute the hyperbolic step:

\begin{corollary}
\label{th:lshs}
Take $f\in\Hol(\D,\D)$ and let $h\in\Hol(\D,\D)$ be a left straightening of $\{f^n\}$. Then for every $z$, $w\in\D$ we have
\begin{equation}
\lim_{n\to+\infty}\omega\bigl(f^n(z),f^n(w)\bigr)=\omega\bigl( h(z),h(w)\bigr)\;.
\label{eq:lshs}
\end{equation}
In particular, $s^f(z)=\omega\bigl(h(z),h(f(z))\bigr)$.
\end{corollary}

\begin{proof}
Let $\{\gamma_n\}\subset\Aut(\D)$ be such that $\gamma_n^{-1}\circ f^n\to h$. Then
\[
\omega\bigl(f^n(z),f^n(w)\bigr)=\omega\bigl(\gamma_n^{-1}\circ f^n(z),\gamma_n^{-1}\circ f^n(w)\bigr)\to \omega\bigl(h(z),h(w)\bigr)
\]
for any $z$, $w\in\D$ and \eqref{eq:lshs} is proved. The final assertion follows immediately by taking $w=f(z)$.
\end{proof}

Armed with this result, we can now answer the question posed above on the hyperbolic step. If $f\in\Hol(\D,\D)$ has a fixed point $z_0\in\D$, the computation of the hyperbolic step is trivial. Indeed, if $f\equiv\id_\D$, then $s^f\equiv 0$. If $f$ is an automorphism, then $s^f(z)=\omega\bigl(z,f(z)\bigr)$ for all $z\in\D$; in particular, $s^f(z)=0$ if and only if $z=z_0$. Finally, if $f$ is not an automorphism, then clearly $s^f\equiv 0$, because all orbits converge to~$z_0$. 

So the only interesting case is when $f$ has no fixed points.

\begin{theorem}
\label{th:zerohypst}
Let $f\in\Hol(\D,\D)$ be without fixed points. Then the following statements are equivalent:
\begin{enumerate}
\item[{\rm (i)}] there exists $z_0\in\D$ such that $s^f(z_0)=0$;
\item[{\rm (ii)}] any left straightening of $\{f^n\}$ is constant;
\item[{\rm (iii)}] $\lim\limits_{n\to+\infty}\omega\bigl(f^n(z),f^n(w)\bigr)=0$ for all $z$, $w\in\D$;
\item[{\rm (iv)}] $s^f\equiv 0$.
\end{enumerate}
\end{theorem}

\begin{proof}
That (ii) implies (iii) follows immediately from \eqref{eq:lshs}. Taking $w=f(z)$ we see that (iii) implies (iv). Moreover, (iv) trivially implies (i); so we are left with proving that (i) implies (ii).

%
Since all left straightenings of~$\{f^n\}$ are constant if and only if one of them is, we can choose $\{\gamma_n\}\subset\Aut(\D)$ such that $\{\gamma_n^{-1}\circ f^n\}$ converges 
to a left straightening $h\in\Hol(\D,\D)$ with $\gamma_n^{-1}\circ f^n(z_0)=h(z_0)=0$; put $H_n=\gamma_n^{-1}\circ f^n$. 

First of all, notice that  for any $z\in\D$ we have
\[
\begin{aligned}
\omega\bigl(0,H_n(z)\bigr)&=\omega\bigl(\gamma_n^{-1}\circ f^n(z_0), \gamma_n^{-1}\circ f^n(z)\bigr)=\omega\bigl(f^n(z_0),f^n(z)\bigr)\\
&\le\omega\bigl(f^{n-1}(z_0),f^{n-1}(z)\bigr)=\omega\bigl(0,H_{n-1}(z)\bigr)\;.
\end{aligned}
\]
So $\{|H_n|\}$ converges non-increasingly to~$|h|$; in particular, $|H_n|\ge|h|$ for all $n\in\N$.

Next, by Corollary~\ref{th:lshs}, $s^f(z_0)=0$ if and only if $\omega\bigl(h(z_0),h\bigl(f(z_0)\bigr)\bigr)=0$; since $h(z_0)=0$, we must have $h\bigl(f(z_0)\bigr)=0$. Notice that $f$ has no fixed points; so $f(z_0)$ is a distinct zero of~$h$. We claim that then $h\equiv 0$. 

Assume, by contradiction, that $h$ is not identically zero. Then $f(z_0)$ is an isolated zero of~$h$; choose a small $r>0$ so that, if $B=B\bigl(f(z_0),r\bigr)\subset\D$ is the Euclidean ball of center~$f(z_0)$ and radius~$r$, then $f(z_0)$ is the unique zero of~$h$ in~$\overline{B}$. Let $\varepsilon=\inf\limits_{\zeta\in\partial B}|h(\zeta)|>0$. Since 
$H_n\to h$ uniformly on compact subsets, we have $\sup\limits_{\zeta\in\partial B}|H_n(\zeta)-h(\zeta)|<\varepsilon$ for $n$ large enough. By Rouch\'e theorem,
$H_n$ and $h$ must eventually have the same number of zeroes in~$B$, counted with multiplicities. But $h$ in~$B$ vanishes only in~$f(z_0)$; since $|H_n|\ge|h|$, the same must happen for~$H_n$ too as soon as $n$ is large enough. 

We have then proved that 
\[
\gamma_n^{-1}\circ f^n\bigl(f(z_0)\bigr)=H_n\bigl(f(z_0)\bigr)=0=H_n(z_0)=\gamma_n^{-1}\circ f^n(z_0)
\] 
for $n$ large enough.
But this implies that $f\bigl(f^n(z_0)\bigr)=f^n(z_0)$, that is, that $f^n(z_0)$ is a fixed point for~$f$, impossible. The contradiction stems from having assumed that $h\not\equiv 0$; so we must have $h\equiv 0$ 
and we are done.
\end{proof}

\section{Parabolic dynamics}
\label{section3}

We now focus on parabolic maps. We recall that a parabolic self-map of~$\D$ is a $f\in\Hol(\D,\D)$ with Wolff point~$\tau_f\in\partial\D$ and such that $f'(\tau_f)=1$.
It turns out that it is easier to work with parabolic (and hyperbolic) self-maps in the setting of the upper half-plane; so we shall first of all restate our hypotheses in this setting.

Let $\mathbb{H}^+=\{w\in\C\mid \mathop{\mathrm{Im}} w>0\}$ be the upper half-plane in the complex plane. 
The \emph{Cayley transform} $\Psi\colon\D\to\mathbb{H}^+$,
given by $\Psi(z)=i\frac{1+z}{1-z}$, is a biholomorphism between $\D$ and $\mathbb{H}^+$ which extends to a homeomorphism between $\overline{\D}$ and the closure
$\overline{\mathbb{H}^+}$ of $\mathbb{H}^+$ in the Riemann sphere by setting $\Psi(1)=\infty$. Using the Cayley transform, we can define the Poincar\'e distance
$\omega_{\mathbb{H}^+}\colon\mathbb{H}^+\times\mathbb{H}^+\to\R^+$ by 
\[
\omega_{\mathbb{H}^+}(w_1,w_2)=\omega\bigl(\Psi^{-1}(w_1),\Psi^{-1}(w_2)\bigr)=\tanh^{-1}\left|\frac{w_2-w_1}{w_2-\overline{w_1}}\right|\;;
\]
we clearly recover the Schwarz-Pick theorem for holomorphic self-maps of~$\mathbb{H}^+$. In particular, we can define the hyperbolic step~$s^F$ of a $F\in\Hol(\mathbb{H}^+,\mathbb{H}^+)$ and all the results of Section~\ref{section2} hold in this setting too.

The Cayley transform induces a bijection between $\Hol(\D,\D)$ and $\Hol(\mathbb{H}^+,\mathbb{H}^+)$: given a map $F\in\Hol(\mathbb{H}^+,\mathbb{H}^+)$ then $f=\Psi^{-1}\circ F\circ\Psi\in\Hol(\D,\D)$ and conversely. We can then say that $F\in\Hol(\mathbb{H}^+,\mathbb{H}^+)$ is parabolic (or hyperbolic) with Wolff point $\tau_F\in\partial\mathbb{H}^+=\R\cup\{\infty\}$ if and only if $f=\Psi^{-1}\circ F\circ\Psi$ is parabolic (or hyperbolic) with Wolff point $\Psi^{-1}(\tau_F)\in\partial\D$.

Combining the Wolff-Denjoy Theorem~\ref{th:WD} with the Julia-Wolff-Carath\'eodory theorem for the upper half-plane (see, e.g., \cite[Corollary 2.3.4]{Ab2022}) we obtain the following statement:

\begin{theorem}
\label{th:WDH}
Let $F\in\Hol(\mathbb{H}^+,\mathbb{H}^+)$ be a holomorphic self-map without fixed points. Then there exists a point $\tau_F\in\partial\mathbb{H}^+$ such that $F^n\to\tau_F$ uniformly on compact subsets. Furthermore, if $\tau_F=\infty$ then 
there exists a real number $F'(\infty)\ge 1$ such that for any sequence $\{w_n\}$ converging non-tangentially to~$\infty$ we have
\begin{equation}
\lim_{n\to\infty}\frac{F(w_n)}{w_n}=\lim_{n\to\infty}F'(w_n)=F'(\infty)\;.
\label{eq:JWC}
\end{equation}
\end{theorem}

Here, a sequence $\{w_n\}\subset\mathbb{H}^+$ converges non-tangentially to~$\infty$ if and only if there exists $\varepsilon>0$ such that $\Im w_n\ge\varepsilon|w_n|$ 
for all $n\in\N$ or, equivalently, if and only if there exists $\delta>0$ such that $\arg w_n\in[\delta,\pi-\delta]$ for all~$n\in\N$; see, e.g., \cite[Proposition 2.2.7]{Ab2022}.

Our strategy for studying the dynamics of a parabolic map $f\in\Hol(\D,\D)$ will then be the following. First of all, up to conjugating $f$ by a rotation, that is up to replacing $f$
by $f_1(z)=\tau_f^{-1}f(\tau_f z)$, we can assume that $\tau_f=1$. Then $F=\Psi\circ f\circ\Psi^{-1}$ will be parabolic with Wolff point at~$\infty$. We shall then study the dynamics  of~$F$ and we shall finally translate the results back to~$f$.

The first useful result that we want to prove is the following lemma, which is a generalization (and with a proof not depending on models) of \cite[Corollary 4.6.10]{Ab2022}.

\begin{lemma}
\label{th:Val0}
Let $F\in\Hol(\mathbb{H}^+,\mathbb{H}^+)$ be parabolic with Wolff point at~$\infty$. Assume that there are a $w_0\in\mathbb{H}^+$ and a subsequence $\{F^{n_k}\}$ 
such that $\{F^{n_k}(w_0)\}$ converges non-tangentially to~$\infty$. Then~$s^F\equiv 0$.
\end{lemma}

\begin{proof}
Put $w_n=F^n(w_0)$ and write $p(w)=F(w)-w$. First of all, we have
\begin{equation}
\tanh\omega_{\mathbb{H}^+}(w_n,w_{n+1})=\left|\frac{w_{n+1}-w_n}{w_{n+1}-\overline{w_n}}\right|=\left|\frac{p(w_n)}{2i\Im w_n+p(w_n)}\right|\le\frac{|p(w_n)/w_n|}{2\frac{\Im w_n}{|w_n|}-\left|\frac{p(w_n)}{w_n}\right|}\;.
\label{eq:par1}
\end{equation}
Since $\{w_{n_k}\}$ converges non-tangentially to~$\infty$, we can find $\varepsilon>0$ such that $\Im w_{n_k}\ge\varepsilon |w_{n_k}|$ for all $k\in\N$. So \eqref{eq:par1} yields
\begin{equation}
\tanh\omega_{\mathbb{H}^+}(w_{n_k},w_{n_k+1})\le\frac{|p(w_{n_k})/w_{n_k}|}{2\varepsilon-\left|\frac{p(w_{n_k})}{w_{n_k}}\right|}\;.
\label{eq:par2}
\end{equation}
Moreover, \eqref{eq:JWC} yields $F(w_{n_k})/w_{n_k}\to 1$;
hence $p(w_{n_k})/w_{n_k}\to 0$ and thus $\omega_{\mathbb{H}^+}(w_{n_k},w_{n_k+1})\to 0$, by~\eqref{eq:par2}. But we have already remarked that $\{\omega_{\mathbb{H}^+}(w_{n},w_{n+1})\}$ is a non-increasing sequence converging to~$s^F(w_0)$; therefore we must have $s^F(w_0)=0$ and, thus, $s^F\equiv 0$ by Theorem~\ref{th:zerohypst}.
\end{proof}

We can now prove the following 
statement, that contains the answer to the question posed by Valiron \cite[p. 148]{Val1954} mentioned in the introduction in the unit disk setting. 

\begin{theorem}
\label{th:Val1}
Let $F\in\Hol(\mathbb{H}^+,\mathbb{H}^+)$ be parabolic with Wolff point at~$\infty$. Fix $w_0\in\mathbb{H}^+$. Then
\[
\frac{F^n}{F^n(w_0)}\to 1
\]
uniformly on compact subsets.
\end{theorem}

\begin{proof}
%
%
Put $w_n=F^n(w_0)$. 
Assume, by contradiction, that $\{w_n^{-1}F^n\}$ does not converge to~1. Then Valiron (see~\cite[pp. 146-148]{Val1954} and Remark~\ref{rem:Val} below) found a subsequence 
$\{n_k\}$ such that $\{w_{n_k}\}$ converges non-tangentially to~$\infty$ and such that $\{|w_{n_k}|^{-1} F^{n_k}\}$ converges to a non-constant map $\psi\in\Hol(\mathbb{H}^+,\mathbb{H}^+)$ satisfying $\psi\circ F=\psi$. 

By Lemma~\ref{th:Val0} and Theorem~\ref{th:zerohypst}, all left straightenings of~$\{F^n\}$ must be constant. Let $\{\gamma_n\}\subset\Aut(\mathbb{H}^+)$ be a sequence such that $\gamma_n^{-1}\circ F^n$ converges to a left straightening~$h$ of~$\{F^n\}$. 
Choose $z$,~$w\in\mathbb{H}^+$ with $\psi(z)\ne \psi(w)$. Observe that
\[
\begin{aligned}
\omega_{\mathbb{H}^+}\bigl(\gamma_n^{-1}\circ F^n(z),\gamma_n^{-1}\circ F^n(w)\bigr)
&=\omega_{\mathbb{H}^+}\bigl(F^n(z),F^n(w)\bigr)\\
&\ge \omega_{\mathbb{H}^+}\bigl(\psi\bigl(F^n(z)\bigr),\psi\bigl(F^n(w)\bigr)\bigr)
=\omega_{\mathbb{H}^+}\bigl(\psi(z),\psi(w)\bigr)>0\;.
\end{aligned}
\]
Now, $\omega_{\mathbb{H}^+}\bigl(\gamma_n^{-1}\circ F^n(z),\gamma_n^{-1}\circ F^n(w)\bigr)\to \omega_{\mathbb{H}^+}\bigl(h(z),h(w)\bigr)$. Hence $h(z)\ne h(w)$ and $h$ is not constant, contradiction.
\end{proof}

\begin{remark}
\label{rem:Val0}
The argument used in the proof of the previous theorem is a particular case of \cite[Corollary~4.6]{AbSh2025}.
\end{remark}

\begin{remark}
\label{rem:Val}
For the sake of completeness, we describe here how Valiron in~\cite[pp. 146-148]{Val1954} produced the subsequence $\{n_k\}$ and the non-constant map~$\psi$ used in the proof of Theorem~\ref{th:Val1}.

First of all, the sequence $\{w_n^{-1}F^n\}$ is normal. Indeed, from every subsequence we can extract a further subsequence $\{w_{n_k}^{-1}F^{n_k}\}$ such that $\arg w_{n_k}$ has a limit in~$[0,\pi]$; therefore the image of $w_{n_k}^{-1}F^{n_k}$ eventually avoids a fixed sector and, thus, $\{w_{n_k}^{-1}F^{n_k}\}$ admits a converging subsequence.

Since $w_n^{-1}F^n(w_0)=1$ for all $n\in\N$, if $\psi_0$ is a constant limit of a subsequence of $\{w_n^{-1}F^n\}$, then $\psi_0\equiv 1$. Assume, by contradiction, that $\{w_n^{-1}F^n\}$ does not converge to the constant~1. Then, by normality, there must exists a subsequence $\{w_{n_k}^{-1}F^{n_k}\}$ converging to a non-constant map~$\psi_1$.
Up to a subsequence, we can also assume that $\arg w_{n_k}\to\phi_\infty\in[0,\pi]$; then
\[
\frac{F^{n_k}(w)}{|w_{n_k}|}=\frac{F^{n_k}(w)}{w_{n_k}}\frac{w_{n_k}}{|w_{n_k}|}\to e^{i\phi_\infty}\psi_1(w)\;.
\]
This means that the sequence $\{|w_{n_k}|^{-1}F^{n_k}\}$ converges to a non-constant holomorphic self-map $\psi=e^{i\phi_\infty}\psi_1$ of~$\mathbb{H}^+$. As a consequence, for any $w\in\mathbb{H}^+$ we have 
\[
\lim_{k\to+\infty}\arg F^{n_k}(w)=\lim_{k\to+\infty}\arg\frac{F^{n_k}(w)}{|w_{n_k}|}=\arg \psi(w)\in(0,\pi)\;.
\]
In other words, $\{F^{n_k}(w)\}$ converges to $\infty$ non-tangentially for all $w\in\mathbb{H}^+$. By \eqref{eq:JWC}, we then have
\[
\lim_{k\to+\infty}\frac{F\bigl(F^{n_k}(w)\bigr)}{F^{n_k}(w)}=1
\]
for all $w\in\mathbb{H}^+$. It follows that 
\[
\psi\bigl(F(w)\bigr)=\lim_{k\to+\infty}\frac{F^{n_k}\bigl(F(w)\bigr)}{|w_{n_k}|}=\lim_{k\to+\infty}\frac{F\bigl(F^{n_k}(w)\bigr)}{F^{n_k}(w)}\cdot\frac{F^{n_k}(w)}{|w_{n_k}|}=\psi(w)\;,
\]
that is, $\psi\circ F\equiv \psi$ and we are done.
\end{remark}

If we translate this result back to the unit disk we obtain the answer promised in the introduction.

\begin{corollary}
\label{th:Val1D}
Let $f\in\Hol(\D,\D)$ be parabolic with Wolff point~$\tau_f$. Fix $z_0\in\D$. Then
\[
\frac{f^n-\tau_f}{f^n(z_0)-\tau_f}\to 1
\]
uniformly on compact subsets.
\end{corollary}

\begin{proof}
Up to conjugating by a rotation, without loss of generality we can assume $\tau_f=1$. Let $F=\Psi\circ f\circ\Psi^{-1}$, where $\Psi$ is the Cayley transform. Then $F$ is a parabolic self-map of $\mathbb{H}^+$ with Wolff point at~$\infty$. Recalling that $\Psi^{-1}(w)=\frac{w-i}{w+i}$, for any $z\in\D$ we find
\[
f^n(z)-1=\frac{F^n\bigl(\Psi(z)\bigr)-i}{F^n\bigl(\Psi(z)\bigr)+i}-1=\frac{-2i}{F^n\bigl(\Psi(z)\bigr)+i}
\]
and hence
\[
\frac{f^n(z)-1}{f^n(z_0)-1}=\frac{F^n(w_0)+i}{F^n\bigl(\Psi(z)\bigr)+i}=\frac{F^n(w_0)}{F^n\bigl(\Psi(z)\bigr)}\frac{1+i/F^n(w_0)}{1+i/F^n\bigl(\Psi(z)\bigr)}\;,
\]
where $w_0=\Psi(z_0)$.
So the assertion follows from Theorem~\ref{th:Val1}.
\end{proof}

As a consequence, we can immediately prove that if a parabolic map has an orbit converging to the Wolff point along a precise slope, then all orbits converge to the Wolff point 
with the same slope.

\begin{corollary}
\label{th:Val3D}
Let $f\in\Hol(\D,\D)$ be parabolic with Wolff point $\tau_f\in\partial\D$. Assume there is $z_0\in\D$ such that $\frac{f^n(z_0)-\tau_f}{|f^n(z_0)-\tau_f|}\to\sigma\in\partial\D$ as $n\to+\infty$. Then $\frac{f^n(z)-\tau_f}{|f^n(z)-\tau_f|}\to\sigma$ for all $z\in\D$.
\end{corollary}

\begin{proof}
Take $z\in\D$. By Corollary~\ref{th:Val1D} we know that $\frac{f^n(z)-\tau_f}{f^n(z_0)-\tau_f}\to 1$; therefore,
\[
\frac{f^n(z)-\tau_f}{|f^n(z)-\tau_f|}=\frac{f^n(z)-\tau_f}{f^n(z_0)-\tau_f}\left|\frac{f^n(z_0)-\tau_f}{f^n(z)-\tau_f}\right|\frac{f^n(z_0)-\tau_f}{|f^n(z_0)-\tau_f|}\to\sigma\;,
\]
as claimed.
\end{proof}

\begin{corollary}
\label{th:Val3H}
Let $F\in\Hol(\mathbb{H}^+,\mathbb{H}^+)$ be parabolic with Wolff point at~$\infty$. Assume there is $w_0\in\mathbb{H}^+$ such that $\arg F^n(w_0)\to\phi\in[0,\pi]$ as $n\to+\infty$. Then $\arg F^n(w)\to\phi$ for all $w\in\mathbb{H}^+$.
\end{corollary}

\begin{proof}
Notice that, if $w=re^{i\theta}\in\mathbb{H}^+$ and $z=\Psi^{-1}(w)\in\D$, then
\[
\frac{1-z}{|1-z|}=\frac{r^{-1}+ie^{-i\theta}}{|r^{-1}+ie^{-i\theta}|}\;.
\]
In particular, $w$ tends to~$\infty$ in such a way that $\arg w$ converges to~$\phi\in[0,\pi]$ if and only if $z\to 1$ in such a way that $\frac{z-1}{|z-1|}\to -ie^{-i\phi}$.
Then the assertion follows from Corollary~\ref{th:Val3D} applied to $f=\Psi^{-1}\circ F\circ\Psi$, which is parabolic with Wolff point~$1$.
\end{proof}

This result does not imply that the orbits of a parabolic map always converge to the Wolff point along a given slope. Indeed, there are examples of parabolic maps with zero hyperbolic step whose orbits do not converge to the Wolff point tangentially to some direction; see, e.g., \cite[Example 17.5.4]{Braccietalbook}.

On the other hand, if $f$ is parabolic with positive hyperbolic step, by Lemma~\ref{th:Val0} 
no orbit can have a subsequence converging non-tangentially to the Wolff point. However, this does not immediately imply that all orbits converge to the Wolff point with the same slope, because in~$\tau_f$ there are two tangential rays and, in principle, an orbit might jump from one ray to the other. Luckily, this does not happen. This result has already been proved by Pommerenke \cite{Pomm1979}; but using the ideas presented so far we can give an easier proof.

\begin{theorem}
\label{th:Val2}
Let $F\in\Hol(\mathbb{H}^+,\mathbb{H}^+)$ be parabolic with Wolff point at~$\infty$ and positive hyperbolic step. Then either $\arg F^n\to 0$ or $\arg F^n\to\pi$, uniformly on compact subsets.
\end{theorem}

\begin{proof}
We first prove that
\begin{equation}
\lim_{n\to+\infty}\left|F^n(w)-\pi/2\right|=\pi/2
\label{eq:limpi}
\end{equation}
for all $w\in\mathbb{H}^+$. 
If this would not be true, there would exist $w_0\in\mathbb{H}^+$ and a subsequence $\{F^{n_k}\}$ such that $\sup_k|\arg F^{n_k}(w_0)-\pi/2|<\pi/2$. But then $F^{n_k}(w_0)\to\infty$ non-tangentially and then, by Lemma~\ref{th:Val0}, we would have $s^F\equiv 0$, impossible.

Now, fix $w_0\in\mathbb{H}^+$. We claim that either $\arg F^n(w_0)\to 0$ or $\arg F^n(w_0)\to\pi$. If this were not true, by \eqref{eq:limpi} we could find subsequences $\{F^{\mu_k}\}$ and $\{F^{\nu_k}\}$ such that $\arg F^{\mu_k}(w_0)\to 0$ and $\arg F^{\nu_k}(w_0)\to\pi$; up to taking further subsequences, we can also assume $\mu_k<\nu_k<\mu_{k+1}$ for all~$k\in\N$. 
We can now define a continuous curve $\sigma\colon[0,+\infty)\to\mathbb{H}^+$ as follows:
\[
\sigma(t)=\begin{cases}
(1-t)w_0+t F(w_0)&\hbox{if $0\le t\le 1$};\\
F^{n+1}\bigl(\sigma(t-(n+1))\bigr)&\hbox{if $n+1\le t\le n+2$ and $n\in\N$}.
\end{cases}
\]
Then $\arg\sigma$ is a continuous function defined on the connected set~$[0,+\infty)$ and whose image contains both values converging to~$0$ and values converging to~$\pi$; therefore we can find a sequence $\{t_k\}\subset[0,1]$ and a subsequence $\{F^{n_k}\}$ such that $\arg F^{n_k}\bigl(\sigma(t_k)\bigr)=\pi/2$ for all~$k\in\N$.
Up to taking further subsequences, we can assume that $t_k\to\tau\in[0,1]$ and that $\arg F^{n_k}$ converges uniformly on~$\sigma([0,1])$. But then we get
$\lim\limits_{k\to+\infty}\arg F^{n_k}(\tau)=\pi/2$, against \eqref{eq:limpi}.

So for each $w_0\in\mathbb{H}^+$ we have either $\arg F^n(w_0)\to 0$ or $\arg F^n(w_0)\to\pi$. But, by Corollary~\ref{th:Val3H}, all points must converge to~$\infty$ with the same argument and we are done.
\end{proof}

We conclude with the translation to the unit disk of this last result. 

\begin{corollary}
\label{th:Val2D}
Let $f\in\Hol(\D,\D)$ be parabolic with Wolff point $\tau_f\in\partial\D$ and positive hyperbolic step. Then all orbits of $f$ converges to~$\tau_f$ tangentially. More precisely, we have
\[
\lim_{n\to+\infty}\frac{f^n(z)-\tau_f}{|f^n(z)-\tau_f|}=\pm i\tau_f\;,
\]
where the sign is the same for all $z\in\D$.
\end{corollary}

\begin{proof}
Up to conjugating by a rotation, we can without loss of generality assume that $\tau_f=1$. The assertion then follows by Theorem~\ref{th:Val2} arguing as in the proof of Corollary~\ref{th:Val3H}, reversing the roles of $f$ and~$F$.
\end{proof}

\backmatter

\end{document}